\documentclass[10pt,a4paper,twoside]{article}

\usepackage{amssymb,amsmath,amsfonts,amsopn,color,amsthm,enumitem}
\usepackage{graphicx}
\usepackage[colorlinks=true]{hyperref}
\hypersetup{urlcolor=blue, citecolor=red, linkcolor=blue}
\usepackage{mathrsfs}

\newcommand{\R}{\mathbb{R}}

\newcommand{\RT}{{\mathbb{R}^3}}

\renewcommand{\le}{\leqslant}
\renewcommand{\ge}{\geqslant}
\renewcommand{\a }{\alpha }

\renewcommand{\b }{\beta }

\newcommand{\n }{\nabla }
\newcommand{\s }{\sigma }
\renewcommand{\t}{\theta}

\renewcommand{\S}{\Sigma}

\renewcommand{\H}{H^1(\RT)}

\renewcommand{\S}{\mathcal{S}}

\newcommand{\irt }{\int_{\RT}}

\def\bbm[#1]{\mbox{\boldmath $#1$}}

\newtheorem{thm}{Theorem}[section]
\newtheorem{prop}[thm]{Proposition}

\newtheorem{lem}[thm]{Lemma}

\usepackage{latexsym}

\newcommand{\beq}{\begin{equation}}
\newcommand{\eeq}{\end{equation}}
\newcommand{\bdoc}{\begin{document}}
\newcommand{\edoc}{\end{document}}
\newcommand{\be}{\begin{enumerate}}
\newcommand{\ee}{\end{enumerate}}
\newcommand{\bdescr}{\begin{description}}
\newcommand{\edescr}{\end{description}}
\newcommand{\ba}{\begin{array}}
\newcommand{\ea}{\end{array}}

\begin{document}

\title{{\bf Infinitely many positive solutions for a Schr\"{o}dinger-Poisson system
\footnote{P. d'Avenia and A. Pomponio are supported by M.I.U.R. -
P.R.I.N. ``Metodi variazionali e topologici nello studio di
fenomeni non lineari''. G. Vaira is supported by M.I.U.R. -
P.R.I.N. ``Variational Methods and Nonlinear Differential Equations''.}}}

\pagestyle{myheadings} \markboth{{\it P. d'Avenia, A. Pomponio, G. Vaira}} {{\it Nonlinear Schr\"odinger-Poisson system}}

\author{P. d'Avenia\thanks{Dipartimento di Matematica, Politecnico di
Bari, Via E. Orabona 4, I-70125 Bari, Italy, e-mail: {\tt
p.davenia@poliba.it}},
A. Pomponio\thanks{Dipartimento di Matematica, Politecnico di
Bari, Via E. Orabona 4, I-70125 Bari, Italy, e-mail: {\tt
a.pomponio@poliba.it}},
G. Vaira\thanks{S.I.S.S.A. Via Beirut, 2-4, I-34014 Trieste, Italy, e-mail: {\tt vaira@sissa.it}}}

%
\date{}

\maketitle

\begin{abstract}
\noindent We are interested in the existence of infinitely many positive solutions of the Schr\"{o}dinger-Poisson system 
\begin{equation*}
\left\{
  \begin{array}{ll}
  -\Delta u+u+V(|x|)\phi u=|u|^{p-1}u,& x \in {\mathbb R}^{3},\\
  -\Delta \phi=V(|x|)u^{2},& x \in {\mathbb R}^{3},
  \end{array}
\right.
\end{equation*}
where $V(|x|)$ is a positive bounded function, $1<p<5$ and $V(r)$, $r=|x|$ has the following decay property: $$V(r)=\frac{a}{r^m}+O\left(\frac{1}{r^{m+\theta}}\right)$$ with $a>0$, $m>\frac{3}{2}$, $\theta>0$. The solutions founded are non-radial.
 \end{abstract} 

\section{Introduction and main results}\label{introduzione}
\def\theequation{1.\arabic{equation}}\makeatother
\setcounter{equation}{0}

In this paper we consider the following nonlinear Schr\"{o}dinger-Poisson system
\begin{equation}
\left\{
  \begin{array}{ll}
  -\Delta u+u+V(x)\phi u=|u|^{p-1}u,& x \in {\mathbb R}^{3},\\
  -\Delta \phi=V(x)u^{2},& x \in {\mathbb R}^{3},
  \end{array}
\right.
\tag{${\cal SP}$}
\label{eq:SP}
\end{equation}
where $p\in(1, 5)$ and $V:\mathbb R^3\rightarrow \mathbb R$ is a non-negative bounded function.

This kind of problem has been introduced in \cite{BenciFortunato1} and arises in an interesting physical context.
In fact, according to a classical model, the interaction of a charge particle
with an electromagnetic field can be described by coupling the nonlinear
Schr\"odinger and the Maxwell equations and so it is also known as Schr\"odinger-Maxwell system. In our case,  $V(x)$ can be interpreted as a changing pointwise charge distribution. There is a wide literature concerning this type of problem in different situations, see for example \cite{A,ADP,AP,CeramiVaira,DM1,DW,DA,K1,K2,IanniVaira,Ru,RuizVaira,WZ,ZZ}.\\
The equations in \eqref{eq:SP} are the Eulero-Lagrange equations of the $C^2$-functional $G:H^1(\R^3)\times D^{1, 2}(\R^3) \rightarrow \R$ where
\begin{align*}
G(u,\phi)=& \,\frac{1}{2}\irt \left(|\nabla u |^2 + |u|^2\right)\, dx + \frac{1}{2}\irt V(x)\phi u^2\, dx- \frac{1}{4}\irt |\nabla\phi|^2\, dx\\
&\; - \frac{1}{p+1}\irt |u|^{p+1}\, dx,
\end{align*}
and the critical points of $G$ are the solutions of \eqref{eq:SP}. Let us observe that the functional $G$ is unbounded both from above and from below, also modulo compact perturbations. As we shall see in Section 2, for all $u\in H^1(\R^3)$, the Poisson equation in \eqref{eq:SP} admits a unique solution $\phi_u\in D^{1, 2}(\R^3)$. Hence we can reduce ourselves to the study of the $C^2$ one variable functional $I:H^1(\R^3)\rightarrow \R$, defined by
\begin{align*}
I(u)
&=G(u, \phi_u)
\\
&=\frac{1}{2}\irt \left(|\nabla u |^2 + |u|^2\right)\, dx +\frac{1}{4}\int_{\mathbb R^3}V(x)\phi_u u^2\, dx-\frac{1}{p+1}\int_{\mathbb R^3}|u|^{p+1}\, dx.
\end{align*}
The critical points of $I$ are the solutions of the problem 
\begin{equation}
-\Delta u+ u+ V(x)\phi_u u=|u|^{p-1}u.
\label{eq:SP'}
\tag{${\cal SP'}$}
\end{equation}
If $u\in H^1(\mathbb R^3)$ is a solution of \eqref{eq:SP'}, then 
$(u, \phi_u)\in H^1(\mathbb R^3)\times D^{1, 2}(\mathbb R^3)$ is a solution of \eqref{eq:SP} and so we will look for solutions of \eqref{eq:SP'}.

In this paper we assume that $V$ is a radial function, that is $V(x)=V(|x|)=V(r)$.
Moreover we assume that $V$ satisfies the following condition:
\begin{enumerate}[label=(\Alph*), ref=\Alph*]
\setcounter{enumi}{21}
	\item \label{Vipo} there are constants $a>0$, $m>\frac32$, $\theta>0$ such that 
\begin{equation*}
V(r)=\frac{a}{r^m}+O\left(\frac{1}{r^{m+\theta}}\right),
\end{equation*}
as $r\rightarrow +\infty$.
\end{enumerate}
In what follows, without any loss of generality, we assume $a=1$.

The main result of this paper can be stated as follows:

\begin{thm}\label{principale}
If $V$ satisfies (\ref{Vipo}), then the problem \eqref{eq:SP'} has infinitely many non-radial positive solutions.
\end{thm}
To prove Theorem \ref{principale} we will construct solutions with large number of bumps near infinity. In fact, since $V(r)\rightarrow 0$ as $r\rightarrow +\infty$, the solutions of \eqref{eq:SP'} can be approximated by using the solution $U$ of the following limit problem
\begin{equation}\label{P0}
\left\{
\begin{array}{ll}
-\Delta u+u=u^p,& \hbox{in } \mathbb R^3,
\\
u>0, &\hbox{in } \mathbb R^3,
\\
u(x)\rightarrow 0, & \hbox{as } |x|\to +\infty.
\end{array}
\right.
\end{equation}

For any positive integer $k$, let us define
$$
P_j=\left(r\cos\frac{2(j-1)\pi}{k}, r\sin\frac{2(j-1)\pi}{k}, 0\right)\in \mathbb R^3,\qquad j=1,\ldots,k,
$$
with $r\in \left[r_0 k\log k, r_1 k\log k \right]$ for some $r_1>r_0>0$ and
$$
z_{r}(x)=\sum_{j=1}^k U_{P_j}(x),
$$
where $U_{P_j}(\cdot):=U(\cdot - P_j)$.

If $x=(x_1, x_2, x_3)\in \RT$, we set
\[
H_s=\left\{u\in\H\;\vline\; 
\begin{array}{l}
	u \textrm{ is even in } x_2, x_3;\\
	u(r\cos\t, r\sin\t, x_3)=u\left(r\cos\left(\t+ \frac{2\pi j}{k}\right),r\sin\left(\t+ \frac{2\pi j}{k}\right), x_3\right)\\
	j=1,\ldots,k-1
\end{array}
\right\}.
\]

Theorem \ref{principale} is a direct consequence of the following result.
\begin{thm}\label{principale2}
If $V$ satisfies (\ref{Vipo}), then there exists an integer $k_0>0$ such that for all $k\ge k_0$, \eqref{eq:SP'} has a positive solution $u_k$ of the form
$$u_k=z_{r_k}+w_k,$$
where $r_k\in \left[r_0 k\log k, r_1 k \log k \right]$, $w_k\in H_s$ and, as $k\rightarrow +\infty$, $\|w_k\|\rightarrow 0$.
\end{thm}
The proof of Theorem \ref{principale2} relies on a Lyapunov-Schmidt reduction. 
This technique is almost standard in the perturbation methods, prevalently in the presence of a small parameter (see \cite{AmbrMalc}) and it has been applied to Schr\"odinger-Poisson type system by several authors (see, for example, \cite{DW,IanniVaira,RuizVaira}). Here we use $k$ as parameter, that is we use the number of bumps of the solutions in the construction of spike solutions for \eqref{eq:SP'}.
This idea has been introduced by J.~Wei and S.~Yan (\cite{WeiYan}). They consider the nonlinear Schr\"{o}dinger equation $$-\Delta u+V(|x|)u=u^p$$ with $V(r)\rightarrow V_0>0$ as $r\rightarrow +\infty$ and they prove the existence of infinitely many positive non-radial solutions for such equation by using the technique outlined above. This method has also been applied for the study of different problems (see for example \cite{WWY,WY2}). 
In our case, however, many technical difficulties arise due to the presence of the non-local term $\phi_u$ and a more careful analysis of the interaction between the bumps is required.

\section{Notation and preliminaries}\label{preliminari}
\def\theequation{2.\arabic{equation}}\makeatother
\setcounter{equation}{0}
Hereafter we use the following notations:
\begin{itemize}
\item $H^{1}({\mathbb R}^{3})$ is the usual Sobolev space endowed with the standard scalar product and norm
\begin{displaymath}
(u, v)=\int_{\mathbb R^3}[\nabla u \nabla v+uv]dx,\qquad \|u\|^{2}=\int_{{\mathbb R}^3}\left[|\nabla u|^{2}+u^{2}\right]dx;
\end{displaymath}
\item $D^{1,2}({\mathbb R}^{3})$ is the completion of $C_0^{\infty}(\mathbb R^3)$ with respect to the norm
 \begin{displaymath}
 \|u\|^{2}_{D^{1, 2}}=\int_{{\mathbb R}^3}
|\nabla u|^{2}dx;
\end{displaymath}
\item $L^q(\RT)$, $1\le q \le +\infty$ denotes the usual Lebesgue space with the standard norm  $|u|_{q}$;
\item for any $\rho>0$ and for any $z\in \mathbb R^3$, $B_{\rho}(z)$ denotes the ball of radius $\rho$ centered at $z$ and $B_\rho=B_{\rho}(0)$;
\item $C, C', C_i, \bar C$ are various positive constants which may also vary from line to line;
\item $A=O(\frac{1}{h})$ means that there exists a constant $C>0$ such that $|A|\le C\cdot \frac{1}{h}$.
\end{itemize}

Let us summarize some properties of $\phi_u$.
\begin{lem}\label{propphi}
For every $u\in \H$, there
exists a unique $\phi_u\in D^{1,2}(\RT)$ solution of $-\Delta \phi=V(x)u^2$. Such a solution $\phi_u$ is non-negative and the following representation formula holds
\begin{equation*}
\phi_u(x)=\int_{\R^3}\frac{1}{|x-y|}V(y)u^2(y)\, dy.
\end{equation*}
Moreover:
\begin{enumerate}[label=\roman*), ref=\roman*)]
\item \label{lemmaphi1} there exist $C,C'>0$ independent of $u\in\H$ such that
$$\|\phi_u\|_{D^{1,2}}\le C  |u|^2_\frac{12}{5},$$
and
\begin{equation*}
\irt V(x)\phi_u u^2\le C' |u|^4_\frac{12}{5};
\end{equation*}
\item \label{lemmaphi2} if $u\in H_s$, then $\phi_u\in D_s$, where
\[
D_s=\left\{\phi\in D^{1,2}(\RT)\;\vline\; 
\begin{array}{l}
	\phi \textrm{ is even in } x_2, x_3;\\
	\phi(r\cos\t, r\sin\t, x_3)=\phi\left(r\cos\left(\t+ \frac{2\pi j}{k}\right),r\sin\left(\t+ \frac{2\pi j}{k}\right), x_3\right)\\
	j=1,\ldots,k-1
\end{array}
\right\};
\]
\item \label{lemmaphi3} for all $v_1, v_2 \in H^1(\mathbb R^3)$ there holds 
\begin{equation}
\|\phi_{v_1}-\phi_{v_2}\|_{D^{1, 2}}\le C(\|v_1\|+\|v_2\|)\|v_1-v_2\|.
\end{equation}
\end{enumerate}
\end{lem}
\begin{proof}
The existence and uniqueness of $\phi_u \in D^{1, 2}(\R^3)$ is a direct application of the Lax-Milgram theorem. Moreover the inequalities in \ref{lemmaphi1} are well known (see, for instance, \cite{BenciFortunato1,DM1,Ru}) and here we give the proofs only of \ref{lemmaphi2} and \ref{lemmaphi3}.\\
{\it \ref{lemmaphi2}} If $g\in O(3)$ and $u:\RT\to\R$, let us set
\begin{equation}
(T_g u)(x)=u(gx).
\label{eq:actsymm}
\end{equation}
Fixed $g\in O(3)$, if $T_g u=u$, we have that
\[
-\Delta (T_g \phi_u)= T_g(-\Delta\phi_u)=T_g(V(|x|)u^2(x))=V(|gx|)u^2(gx)=V(|x|)u^2(x)=-\Delta \phi_u
\]
and  then $T_g \phi_u=\phi_u$.
\\
Since $H_s$ and $D_s$ are the sets of the fixed points with respect to the action \eqref{eq:actsymm}
for all $g=(g_{i,j})\in O(3)$ with
\[
g_{ij}=
\left\{
\begin{array}{lll}
	-1& &i=j,i=2,3\\
	\delta_{ij}& &\hbox{otherwise}
\end{array}
\right.
\textrm{ or }
g=
\left(
\begin{array}{ccc}
	\cos \frac{2\pi j}{k}&-\sin \frac{2\pi j}{k}&0\\
	\sin \frac{2\pi j}{k}&\cos \frac{2\pi j}{k}&0\\
	0&0&1	
\end{array}
\right),
\]
we conclude that if $u\in H_s$, then $\phi_u\in D_s$.
\\
{\it \ref{lemmaphi3}} We have
\begin{align*}
\|\phi_{v_1}-\phi_{v_2}\|^2_{D^{1, 2}}&= 
\irt |\n (\phi_{v_1}-\phi_{v_2})|^2\,d x
=\irt (v_1^2-v_2^2)(\phi_{v_1}-\phi_{v_2})\, d x
\\
&\le C |v_1^2-v_2^2|_{\frac 65} \|\phi_{v_1}-\phi_{v_2}\|_{D^{1,2}}
\\
&\le C |v_1+v_2|_{\frac{12}5} |v_1-v_2|_{\frac{12}5} \|\phi_{v_1}-\phi_{v_2}\|_{D^{1,2}}
\\
&\le C(\|v_1\|+\|v_2\|)\|v_1-v_2\| \|\phi_{v_1}-\phi_{v_2}\|_{D^{1,2}} .
\end{align*}
\end{proof}

Analogously, the following lemma holds.
\begin{lem}\label{propphi2}
For any $u,v \in H^1(\mathbb R^3)$, let $\bar \phi$ be the solution in $D^{1, 2}(\mathbb R^3)$ of $-\Delta \bar \phi=V(x)uv$. Then 
\begin{equation*}
\|\bar \phi\|_{D^{1, 2}}\le C|V(x)uv|_{\frac 65} \le C \|u\| \|v\|
\end{equation*}
and for any $z,w \in H^1(\mathbb R^3)$
\begin{equation}\label{stimaphi2}
\int_{\R^3}V(x)\bar \phi zw\, dx \le C\|\bar\phi\|_{D^{1, 2}} |V(x)zw|_{\frac 65} \le C \|u\| \|v\| \|z\| \|w\|.
\end{equation}
\end{lem}

Let $\a\in (0,1),$ by (\ref{Vipo}), for any $x\in B_{\a r}$ and for $i=1,\ldots,k$, we have that 
\beq\label{eq:vp}
V(x+P_i)=\frac{1}{|x+P_i|^m}+O\left(\frac{1}{|x+P_i|^{m+\t}}\right).
\eeq
Moreover, for any $x\in B_{\a r}$ and $\b>0$
\beq\label{eq:zp}
\frac{1}{|x+P_i|^\b}=\frac{1}{|P_i|^\b}\left(1+O\left(\frac{|x|}{|P_i|}\right)\right).
\eeq
Let us evaluate now the distance between the various $P_i$, $i=1, \ldots, k$. Indeed, for $i\neq j$, by making a simple computation, we find
\begin{equation}\label{distanzapunti}
|P_i-P_j|=2r\left|\sin\frac{(i-j)\pi}{k}\right|.
\end{equation}

Finally, let us recall the following elementary inequalities which hold for all $a, b, b_1, b_2 \in \R$ and $p>1$:
\begin{equation}
\left|a+b\right|^p\le C\left(\left|a\right|^p+\left|b\right|^p\right)
\label{dis}
\end{equation}
and, if $|a|\le 1$,
\begin{multline} \label{dis1}
\left||a+b_1|^{p}-|a+b_2|^p-p|a|^{p-1}(b_1-b_2)\right| \\
\le 
\left\{
\begin{array}{ll}
C\left(|b_1|^{p-1}+|b_2|^{p-1}\right)|b_1-b_2|,& \hbox{if } p\le 2,\\
C\left(|b_1|^{p-1}+|b_2|^{p-1}+|b_1|+|b_2|\right)|b_1-b_2|,& \hbox{if } p>2
\end{array}
\right.
\end{multline}
where the constant $C$ depends only on $p$.

As mentioned in the introduction, we denote by $U$ the unique positive radial solution in $H^1(\mathbb R^3)$ of the problem (\ref{P0}). This solution satisfies the following decay property (see \cite{Kwong}): $$\lim_{r\rightarrow +\infty}U(r)re^r=C>0,\qquad \lim_{r\rightarrow +\infty}\frac{U'(r)}{U(r)}=-1,\qquad r=|x|,$$ for some constant $C$. \\ The function $U$ is a critical point of the $C^2$ functional $I_0: H^1(\mathbb R^3)\rightarrow \mathbb R$ defined as 
\begin{equation*}
I_0(u)=\frac{1}{2}\|u\|^2-\frac{1}{p+1}\int_{\mathbb R^3}|u|^{p+1}\, dx.
\end{equation*}
Furthermore the solution $U$ is nondegenerate (up to translations). \\More specifically, there holds
\begin{lem}\label{nondegeneratezzaU}
Define the operator $Q:H^1(\mathbb R^3)\rightarrow \mathbb R$ as $$Q[v]:=I_0''(U)[v, v]=\int_{\mathbb R^3}\left[|\nabla v|^2+v^2-pU^{p-1}v^2 \right]\, dx.$$ We denote by $U_j=\frac{\partial U}{\partial x_j}$. Then there hold:
\begin{itemize}
\item $Q[U]=(1-p)\|U\|^2<0$;
\item $Q[U_j]=0$, $j=1, 2, 3$;
\item $Q[v]\ge C\|v\|^2$ for all $v\bot U$, $v\bot U_j$, $j=1, 2, 3$.
\end{itemize}
\end{lem}
For the proof we refer, for instance, to \cite[Lemma 8.6]{AmbrMalc}.

Finally, let us recall the following results (see \cite[Corollary 3.2]{KW} and \cite[Lemma 3.7]{ACR}).
\begin{lem}\label{berlions}
Let $\beta_1\ge 1$ and $\beta_2 \ge 1$ be two positive numbers. Then we have 
\begin{equation*}
\int_{\mathbb R^3}U_{P_i}^{\beta_1}U_{P_j}^{\beta_2}\, dx = O\left(e^{-(\min\{\beta_1, \beta_2\}-\delta)|P_i-P_j|}\right),\qquad i\neq j,
\end{equation*}
where $\delta>0$ is any small number.
\end{lem}

\begin{lem}\label{ACR}
For all $i,j=1,\ldots,k$, $i\neq j$, we have
\[
\irt U_{P_i}^p U_{p_j} \,d x =C \frac{e^{-|P_i-P_j|}}{|P_i-P_j|} (1+o(1)).
\]
\end{lem}

\section{Proof of Theorem \ref{principale2}}

The proof of Theorem \ref{principale2} relies on a Lyapunov-Schmidt reduction. Let
$$
Z_{i, j}=\frac{\partial U_{P_i}}{\partial x_j},\qquad i=1, \ldots , k \hbox{ and } j=1, 2, 3,
$$ 
and define 
$$W:=\left\{w\in H_s\,:\, \int_{\mathbb R^3}U_{P_i}^{p-1}Z_{i, j}w\, dx=0,\,\, i=1,\ldots, k; j=1, 2, 3 \right\}$$
and $P_W$ the orthogonal projection onto $W$. Our approach is to find first a solution $w\in W$ of the auxiliary equation $$P_WI'(z_{r}+w)=0$$ and then to solve the remaining finite dimensional equation, namely the bifurcation equation 
$$
(I-P_W)I'(z_{r}+w)=0.
$$

In what follows we always assume that $V$ satisfies (\ref{Vipo}) and $$r\in S_k:=\left[\left(\frac{m}{\pi}-\beta\right)k\log k,\,\left(\frac{m}{\pi}+\beta\right)k\log k\right],$$ where $m$ is the constant in (\ref{Vipo}) and $\beta>0$ is a small constant.

\subsection{The auxiliary equation}
\def\theequation{3.\arabic{equation}}\makeatother
\setcounter{equation}{0}

In the sequel we find a solution of the auxiliary equation, namely we prove the following proposition.
\begin{prop}\label{ausiliaria}
There exists an integer $k_0>0$, such that for each $k\ge k_0$, there is a $C^1$ map $w: S_k \rightarrow H_s$, $w=w(r)$, satisfying $w\in W$ and $P_W I'(z_r+w)=0$. Moreover, there is a small $\sigma>0$, such that 
\begin{equation*}
\|w\|\le \frac{C}{k^{m-\frac 12+\sigma}}.
\end{equation*}
\end{prop}
We begin with some estimates.

\begin{lem}\label{derivprimaoper}
There exists an integer $k_0>0$, such that for each $k\ge k_0$, there is a small $\sigma>0$ such that
\begin{equation}\label{stimaderprima}
\|I'(z_r)\|\le \frac{C}{k^{m-\frac12+\sigma}}.
\end{equation}
\end{lem}
\begin{proof}
Let $v\in H_s$, with $\|v\|=1$. Taking into account that $U_{P_i}$ are solutions of (\ref{P0}), we have 
$$I'(z_r)[v]=\underbrace{\int_{\mathbb R^3}V(x)\phi_{z_r}z_r v\, dx}_{(I)}-\underbrace{\int_{\mathbb R^3}\left[z_r^p-\sum_{i=1}^k U_{P_i}^p\right]v\, dx}_{(II)}.$$
Let us evaluate separately the two terms. By using H\"{o}lder inequality we obtain:
\begin{align*}
|(I)|&=\left| \sum_{i=1}^k\int_{\mathbb R^3}V(x)\phi_{z_r}U_{P_i}v\, dx\right|\le \sum_{i=1}^k\left[\left(\int_{\mathbb R^3}(V(x)\phi_{z_r}U_{P_i})^2\,dx\right)^{1/2}|v|_2\right]
\\
&\le  C\|\phi_{z_r}\|_{D^{1, 2}}\sum_{i=1}^k\left(\int_{\mathbb R^3}(V(x)U_{P_i})^3\,dx\right)^{1/3}.
\end{align*}
Now, let be $\alpha \in (0, 1)$. By \eqref{eq:vp} and \eqref{eq:zp} and since $U$ decays exponentially outside a ball we find
\begin{align}\label{I.1}
\nonumber
\int_{\mathbb R^3}(V(x)U_{P_i}(x))^3\,dx&=\int_{\mathbb R^3}(V(x+P_i)U(x))^3\,dx\\
\nonumber
&=\int_{B_{\alpha r}}(V(x+P_i)U(x))^3\,dx+\int_{\RT\setminus B_{\alpha r}}(V(x+P_i)U(x))^3\,dx\\
&\le  \frac{C}{r^{3m}},
\end{align}
and 
\begin{align*}
\|\phi_{z_r}\|^2_{D^{1, 2}}&=\int_{\mathbb R^3}|\nabla \phi_{z_r}|^2\, dx=\int_{\mathbb R^3}V(x)\phi_{z_r}z_r^2\, dx
\\
&=\sum_{i=1}^k \int_{\mathbb R^3}V(x)\phi_{z_r}U_{P_i}^2\, dx+\sum_{i\neq j}\int_{\mathbb R^3}V(x)\phi_{z_r}U_{P_i}U_{P_j}\, dx
\\
&\le  C\|\phi_{z_r}\|_{D^{1, 2}}\sum_{i=1}^k \left( \int_{\mathbb R^3}(V(x)U_{P_i}^2)^{6/5}\, dx \right)^{5/6}
\\
&\quad+C\|\phi_{z_r}\|_{D^{1, 2}}\sum_{i\neq j}\left(\int_{\mathbb R^3}(V(x)U_{P_i}U_{P_j})^{6/5}\, dx \right)^{5/6}.
\end{align*}
As done in (\ref{I.1}) we find 
\begin{equation*}
\int_{\mathbb R^3}(V(x)U_{P_i}^2)^{6/5}\, dx \le \frac{C}{r^{6m/5}}.
\end{equation*}
Moreover by Lemma \ref{berlions} and (\ref{distanzapunti}), since $r\in S_k$,  we find
\begin{align*}
\sum_{i\neq j}\left(\int_{\mathbb R^3}(V(x)U_{P_i}U_{P_j})^{6/5}\, dx \right)^{5/6}
&\le  C\sum_{i\neq j} e^{-(1-\delta)|P_i-P_j|}
\\
&=  C k\cdot \sum_{i=2}^k e^{-(1-\delta)|P_i-P_1|}
\\
&\le  Ck\cdot e^{-(1-\delta)\frac{2\pi}{k}r}
\\
&\le C\frac{k}{k^{2(1-\delta)(m-\beta \pi)}}.
\end{align*}
Hence
$$\|\phi_{z_r}\|_{D^{1, 2}}\le  C\cdot \frac{k}{r^m}.$$
Thus, since $r\in S_k$ and $m>\frac32$,
$$|(I)|\le C\cdot \frac{k^2}{r^{2m}}\le \frac{C}{k^{m-\frac12+\sigma}}.$$
Finally we consider (II). These estimates have been done in \cite{WeiYan}; we sketch the proof for the sake of completeness. Let us define for all $j=1,.., k$
$$\Omega_j:=\left\{x=(x', x_3)\in \mathbb R^2 \times \mathbb R \,\, :\,\, \langle \frac{x'}{|x'|},\,\frac{P_j}{|P_j|}\rangle \ge \cos\frac{\pi}{k}\right\}.$$
For any $x\in \Omega_1$:  $U_{P_i}^p\le U_{P_1}^{p-1}U_{P_i}$. Thus, by using \cite[Lemma A.1]{WeiYan}
\begin{align*}
\left|(II)\right|
&= k \left|\int_{\Omega_1}\left(z_r^p-\sum_{i=1}^k U_{P_i}^p\right)v\, dx \right|
\\
&\le  C k\int_{\Omega_1}U_{P_1}^{p-1}\sum_{j=2}^k U_{P_j}|v|\, dx 
\\
&\le  Ck\left( \int_{\Omega_1}U_{P_1}^{\frac{(p-1)(p+1)}{p}}\left(\sum_{j=2}^k U_{P_j}\right)^{\frac{p+1}{p}}\, dx\right)^{\frac{p}{p+1}}\left(\int_{\Omega_1}|v|^{p+1}\, dx\right)^{\frac{1}{p+1}}
\\
&\le  Ck^{\frac{p}{p+1}}\left( \int_{\Omega_1}U_{P_1}^{\frac{(p-1)(p+1)}{p}}\left(\sum_{j=2}^k U_{P_j}\right)^{\frac{1}{p}}\sum_{j=2}^k U_{P_j}\, dx\right)^{\frac{p}{p+1}}
\\
&\le  Ck^{\frac{p}{p+1}}e^{-\frac{\eta r\pi}{(p+1)k}}\left(\sum_{j=2}^k\int_{\mathbb R^3}U_{P_1}^{\frac{p^2-\eta}{p}}U_{P_j}\, dx\right)^{\frac{p}{p+1}}
\\
&\le  C k^{\frac{p}{p+1}}e^{-\frac{\eta r\pi}{(p+1)k}}\left(\sum_{j=2}^k e^{-(1-\delta)|P_1-P_j|}\right)^{\frac{p}{p+1}}
\\
&\le  C k^{\frac{p}{p+1}}e^{-\frac{\pi}{p+1}(\eta+2p(1-\delta))(\frac{m}{\pi}-\beta)\log k}
\\
&\le  C\frac{1}{k^{\left[-\frac{p}{p+1}+\frac{\pi}{p+1}\left(\eta+2p(1-\delta)\right)\left(\frac{m}{\pi}-\beta\right)\right]}}.
\end{align*}
Since $$-\frac{p}{p+1}+\frac{2pm}{p+1}>m-\frac{1}{2},$$ we obtain
$$\left|(II)\right| \le \frac{C}{k^{m-\frac 1 2 +\sigma}}.$$ 
Putting together (I) and (II) we find (\ref{derivprimaoper}).
\end{proof} 
Now we are concerned with the invertibility of the operator $$L:=P_W I''(z_r): W\rightarrow W.$$ 
\begin{lem}\label{inv}
For $k$ sufficiently large $L$ is invertible and $\|L^{-1}\|\le \bar{C}$.
\end{lem}
In order to prove Lemma \ref{inv} let us decompose $W=A\oplus B$ where $$A:=\langle \left\{P_W U_{P_i}\right\}_{i=1...k}\rangle; \qquad B:=\left(A\oplus W^{\bot}\right)^{\bot}.$$  Lemma \ref{inv} follows immediately after showing the next result.
\begin{lem}
For $k$ sufficiently large there exist two positive constants $C_1, C_2$ such that
\begin{enumerate}[label=(\alph*), ref=(\alph*)]
\item $I''(z_r)[u, u] \le -C_1 \|u\|^2$, for all $u\in A$;
\item $I''(z_r)[u, u] \ge C_2\|u\|^2$, for all $u\in B$.
\end{enumerate}
\end{lem}
\begin{proof}
The proof is very similar to \cite[Lemma 3.4]{RuizVaira}. We sketch it for the sake of completeness.\\
Let $u\in A$. Then $$u=\sum_{i=1}^k \lambda_i P_W U_{P_i}, \qquad \lambda_i\in \mathbb R, \quad i=1,\ldots, k.$$
For $i=1,\ldots, k$, $P_W U_{P_i} \in W$. Hence we can write $$P_W U_{P_i}=U_{P_i}-\psi_i, \qquad i=1,\ldots,k,$$ where $\psi_i$ are given by $$\psi_i=\sum_{l,j}(U_{P_i}, Z_{l, j})\frac{Z_{l, j}}{\|Z_{l, j}\|^2},\qquad l\neq i,$$ and the functions $Z_{l, j}$ satisfy $-\Delta Z_{l, j}+Z_{l, j}=pU_{P_l}^{p-1}Z_{l, j}$. 
Since for $k\rightarrow +\infty$ we have $|P_i-P_l|\rightarrow +\infty$, we get $(U_{P_i}, Z_{l, j})=o(1)$ as $k\rightarrow +\infty$. This implies $\|\psi_i\|=o(1)$ as $k\rightarrow +\infty$ for $i=1, ..., k$.\\ Applying the bilinear form given by $I''(z_r)$, using the fact that $I''(z_r)$ maps bounded sets onto bounded sets and that $\|\psi_i\|=o(1)$ we obtain that 
$$
I''(z_r)[u, u]=I''(z_r)\left[\sum_{i=1}^k \lambda_i U_{P_i}, \, \sum_{i=1}^k \lambda_i U_{P_i}\right]+o(1).
$$
Furthermore, by making simple computations, reasoning as in the proof of Lemma \ref{derivprimaoper} and by Lemma \ref{nondegeneratezzaU}, we find
$$
I''(z_r)[u, u]=\sum_{i=1}^k \lambda_i^2I''_0(U_{P_i})[U_{P_i}, U_{P_i}]+o(1)\le (1-p)\sum_{i=1}^k \lambda_i^2\|U_{P_i}\|^2+o(1)<-C_1.
$$ 
Then $I''(z_r)$ is negative definite on $A$.
\\
We now prove that $I''(z_r)$ is positive definite on $B$. Choose an arbitrary $u\in B$ and we assume, for simplicity, that $\|u\|=1$.  
We denote by $\hat{\phi}$ the solution of $-\Delta
\hat{\phi}=V(x)z_{r}u$. Then
$$I''(z_r)[u, u]=\int_{\mathbb R^3}\left[|\nabla u|^2+u^2+V(x)\phi_{z_r}u^2+2V(x)\hat{\phi}z_r u-pz_r^{p-1}u^2\right]\, dx.$$
Since $u$ is bounded and reasoning as in Lemma \ref{derivprimaoper} we find for $\alpha\in (0, 1)$,
\begin{align*}
\int_{\mathbb R^3}V(x)\phi_{z_r}u^2\, dx&= \int_{\mathbb R^3}V(x)\phi_u z_r^2\, dx\\
&= \sum_{i=1}^k \int_{\mathbb R^3} V(x)\phi_u U_{P_i}^2\, dx+\sum_{i\neq j}\int_{\mathbb R^3} V(x)\phi_u U_{P_i}U_{P_j}\, dx\\
& \le C\sum_{i=1}^k\left(\int_{\mathbb R^3} \left(V(x)U_{P_i}^2\right)^{6/5}\, dx\right)^{5/6}+C\sum_{i\neq j}\left(\int_{\mathbb R^3} \left(V(x)U_{P_i}U_{P_j}\right)^{6/5}\, dx\right)^{5/6}\\
&=o(1),
\end{align*}
for $k$ sufficiently large. In the same way one can prove that for $k\rightarrow+\infty$ $$\int_{\mathbb R^3}V(x)\hat{\phi}z_r u\, dx=o(1).$$
As done in Lemma \ref{derivprimaoper}, it can be proved that
$$\int_{\mathbb R^3}z_{r}^{p-1}u^2dx= \int_{\mathbb R^3}\sum_{i=1}^k U_{P_i}^{p-1}u^2dx+o(1).$$ Hence
\[
I''(z_{r})[u, u]= \int_{\mathbb R^3}\left[|\nabla
u|^2+u^2-p\sum_{i=1}^k U_{P_i}^{p-1}u^2\right]\, dx+o(1).
\]
At this point, arguing, for example, as in \cite{RuizVaira}, we have that 
\begin{equation*}
I''(z_{r})[u, u]\ge C_2>0,
\end{equation*}
and (b) follows.
\end{proof}

We are now ready to prove Proposition \ref{ausiliaria}.
\begin{proof}[Proof of Proposition \ref{ausiliaria}]
Let us consider $J(w)=I(z_{r}+w)$, $w\in W$ and expand it as follows:
\begin{align*}
J(w)&=I(z_r)+I'(z_r)[w]+\frac{1}{2}I''(z_r)[w, w]+R_{z_r}(w)\\
&= J(0)+l(w)+\frac{1}{2}\langle L w,\, w \rangle+R_{z_r}(w),
\end{align*}
where 
$$l(w)=I'(z_r)[w], \qquad \langle L w,\, w \rangle=I''(z_r)[w, w],$$ and 
\begin{align*}
R_{z_r}(w)&=\frac{1}{4}\int_{\mathbb R^3}V(x)\phi_w w^2\, dx+\int_{\mathbb R^3}V(x)\phi_w z_r w\, dx
\\
&\quad-\frac{1}{p+1}\int_{\mathbb R^3}\left[|z_r+w|^{p+1}-z_r^{p+1}-\frac{p(p+1)}{2}z_r^{p-1}w^2-(p+1)z_r^pw\right] dx.
\end{align*}
Since $l(w)$ is a bounded linear functional in $W$, by Riesz Theorem there exists an $l_k\in W$ such that $$l(w)=\langle l_k, w \rangle.$$
Now we want to find a critical point of $J$, that is a $w\in W$ such that 
\begin{equation}\label{pfisso}
0=J'(w)=l_k+Lw+R_{z_r}'(w).
\end{equation}
Since by Lemma \ref{inv} $L$ is invertible, we can rewrite (\ref{pfisso}) in the following way $$w=A(w):=-L^{-1}l_k-L^{-1}R'_{z_r}(w).$$ Thus the problem of finding a critical point of $J(w)$ is equivalent to find a fixed point of $A$. To this end, let 
$$
B:=\left\{w\in W:\,\,\, \|w\|\le \frac{C}{k^{m-\frac12+\sigma}}\right\},
$$ 
where $\sigma>0$ is small. We have to prove that $A(B)\subset B$ and that $A$ is a contraction in $B$.\\
Let $w\in B$. By using Lemmas \ref{derivprimaoper} and \ref{inv} 
\begin{align*}
\|A(w)\| & \le \| L^{-1}\|\left(\| l_k \|+\| R'_{z_r}(w) \|\right)\\
& \le \bar{C}\left(\frac{1}{k^{m-\frac{1}{2}+\sigma}}+\| R'_{z_r}(w) \|\right).
\end{align*}
Let us evaluate $\| R'_{z_r}(w) \|$. We denote by $\bar{\phi}$ the solution in $D^{1, 2}(\R^3)$ of $-\Delta \bar{\phi}= V(x)v w$, then 
\begin{align*}
\| R'_{z_r}(w) \| &= \sup_{\|v\|=1}\left| R'_{z_r}(w)[v] \right|\\
&\le\sup_{\|v\|=1}\left[ \int_{\mathbb R^3}V(x)\phi_w |wv|\, dx+2\int_{\mathbb R^3}V(x)|\bar{\phi}z_rw|\, dx+\int_{\mathbb R^3}V(x)\phi_w z_r|v|\, dx\right.\\
& \quad \left.+\int_{\R^3}\left||z_r+w|^p-p z_r^{p-1}w -z_r^p\right|\cdot |v|\, dx \right].
\end{align*}
If $p\le 2$ then, by (\ref{stimaphi2}) and (\ref{dis1}),
\begin{align*}
\| R'_{z_r}(w) \| &\le C\sup_{\|v\|=1}\left[\|w\|^2+\int_{\R^3}|w|^p |v|\, dx\right] \\
& \le C (\|w\|^2+\|w\|^p) \le C\|w\|^p.
\end{align*}
Thus, since $w\in B$,
$$\|A(w)\|\le \widetilde{C}\left(\frac{1}{k^{m-\frac{1}{2}+\sigma}}+\frac{1}{k^{p(m-\frac{1}{2}+\sigma)}}\right)\le  \frac{C}{k^{m-\frac{1}{2}+\sigma}}.$$
If, now, $p>2$ then, by using again (\ref{stimaphi2}) and (\ref{dis1}),
\[
\|R'_{z_r}(w)\|\le  C\sup_{\|v\|=1}\left(\|w\|^2+\int_{\mathbb R^3}\left(|w|^2+|w|^{p}\right)|v|\, dx\right)\le C\|w\|^2
\]
and then 
$$\|A(w)\|\le\frac{C}{k^{m-\frac{1}{2}+\sigma}}.$$
Hence, in both cases, $A$ maps $B$ into $B$.\\
Finally, let us prove that $A$ is a contraction. Let be $w_1$, $w_2 \in B$. Then
\begin{align*}
\| A(w_1)-A(w_2)\|&\le  \|L^{-1}\|(\|l_k\|\|w_1-w_2\|+\|R'_{z_r}(w_1)-R'_{z_r}(w_2)\|)
\\
&\le \bar{C} \left(\frac{C}{k^{m-\frac{1}{2}+\sigma}}\|w_1-w_2\|+\|R'_{z_r}(w_1)-R'_{z_r}(w_2)\|\right).
\end{align*}
Moreover we have
\begin{align*}
\|R'_{z_r}(w_1)-R'_{z_r}(w_2)\|&=\sup_{\|v\|=1}\Big|R'_{z_r}(w_1)[v]-R'_{z_r}(w_2)[v]\Big|
\\
&\le  \sup_{\|v\|=1}\left[ \int_{\mathbb R^3}V(x)\phi_{w_1}|w_1-w_2|\cdot|v|\, dx
+2\irt V(x)|\bar{\phi}_1|z_r|w_1-w_2|\, dx\right.
\\
& \quad +\int_{\mathbb R^3}V(x)|\phi_{w_1}-\phi_{w_2}|z_r|v|\, dx
+2\int_{\mathbb R^3}V(x)|\bar{\phi}_1-\bar{\phi}_2|z_r|w_2|\, dx
\\
&\quad +\int_{\mathbb R^3}V(x)|\phi_{w_1}-\phi_{w_2}|\cdot|w_2 v|\, dx
\\
&\quad \left. +\int_{\mathbb R^3}\left||z_r+w_1|^p-|z_r+w_2|^p-pz_r^{p-1}(w_1-w_2)\right|\cdot |v|\, dx\right],
\end{align*}
where $\bar{\phi}_i$ is the solution in $D^{1, 2}(\mathbb R^3)$ of $-\Delta\bar{\phi}_i=V(x)z_rw_i$, $i=1, 2$. If now $p\le 2$, by \ref{lemmaphi3} of Lemma \ref{propphi} and (\ref{dis1})
\begin{align*}
\| R'_{z_r}(w_1)-R'_{z_r}(w_2)\|
& \le  C \sup_{\|v\|=1}\left[\frac{1}{k^{m-\frac{1}{2}+\sigma}}\|w_1-w_2\|
+\int_{\mathbb R^3}\left(|w_1|^{p-1}+|w_2|^{p-1}\right)|w_1-w_2|\cdot |v|\, dx\right]
\\
&\le C \left[\frac{1}{k^{m-\frac{1}{2}+\sigma}}
+\frac{1}{k^{(p-1)(m-\frac{1}{2}+\sigma)}}\right]\|w_1-w_2\|
\\
&\le  \frac{C}{k^{(p-1)(m-\frac{1}{2}+\sigma)}}\|w_1-w_2\|
\end{align*}
If $p>2$ by using again \ref{lemmaphi3} of Lemma \ref{propphi} and (\ref{dis1})
\begin{align*}
\|R'_{z_r}(w_1)-R'_{z_r}(w_2)\|
&\le  C\sup_{\|v\|=1}\left[ \frac{1}{k^{m-\frac{1}{2}+\sigma}}\|w_1-w_2\|\right.
\\
&\quad \left. +\int_{\mathbb R^3}\left(|w_1|^{p-1}+|w_2|^{p-1}+|w_1|+|w_2|\right)|w_1-w_2|\cdot |v|\, dx\right]\\
&\le  \frac{C}{k^{(m-\frac{1}{2}+\sigma)}}\|w_1-w_2\|.
 \end{align*}
 In both cases $A$ is a contraction since for $k$ large $\frac{C}{k^{s(m-\frac12+\sigma)}}\ll 1$ for $s=1$ and $s=p-1$.
\end{proof}

\subsection{The reduced functional}
This section is devoted to solve the finite dimensional equation, namely the bifurcation equation. \\ To this aim,
let us define the reduced functional $F:S_k \rightarrow \mathbb R$ such that for all $r\in S_k$ $$F(r)=I(z_r+w),$$ where $w=w(r,k)$ is the unique solution of the auxiliary equation. By Proposition \ref{ausiliaria} let us recall that 
$$
\|w\|\le \frac{C}{k^{m-\frac12+\s}}.
$$
By Lemma \ref{derivprimaoper} and since $I''$ maps bounded sets onto bounded sets then 
\begin{align*}
F(r)&=I(z_r)+I'(z_r)[w]+\frac{1}{2}I''(\xi)[w, w]\\
&=I(z_r)+O\left(\|w\|^2\right) \\ 
&= I(z_r)+k\cdot O\left(\frac{1}{k^{2m+\s}}\right),
\end{align*}
where $\sigma>0$ is small.
Then by Proposition \ref{Izr} the reduced functional is given by
\begin{equation*}
F(r)=k\left[C_0+\frac{B_1}{r^{2m}}+ \frac{B_2 k\log k}{r^{2m+1}}-B_3 \sum_{i=2}^k\int_{\mathbb R^3}U_{P_1}^p U_{P_i}\, dx+O\left(\frac{1}{k^{2m+\sigma}}\right)\right],
\end{equation*}
where $C_0, B_1, B_2, B_3$ are positive constants. The problem 
\beq\label{eq:max}
\max\left\{F(r): r\in S_k\right\}
\eeq
has a solution since $F$ is continuous on a compact set. We have to show that this maximum is an interior point of $S_k$.\\
Let us denote with $\bar F$ the function
\[
\bar F (r):=\frac{B_1}{r^{2m}} +\frac{B_2 k\log k}{r^{2m+1}}-B_3 \sum_{i=2}^k\int_{\mathbb R^3}U_{P_1}^p U_{P_i}\, dx.
\]
By Lemmas \ref{ACR} and \ref{A6}
\[
C_1 \frac{e^{- \frac{2\pi r}{k}}}{\log k}\le\sum_{i=2}^k\int_{\mathbb R^3}U_{P_1}^p U_{P_i}\, dx\le C_2 e^{- \frac{2\pi r}{k}}.
\]
So we define
\begin{align*}
F_1(r)&:=\frac{B_1}{r^{2m}} +\frac{B_2 k\log k}{r^{2m+1}}-B_5 e^{- \frac{2\pi r}{k}},\\
F_2(r)&:=\frac{B_1}{r^{2m}} +\frac{B_2 k\log k}{r^{2m+1}}-B_4 \frac{e^{- \frac{2\pi r}{k}}}{\log k}.
\end{align*}
For $k$ sufficiently large,
\[
F_1(r) \le \bar F(r) \le F_2(r)\quad\hbox{ in }S_k
\]
and, moreover, we have
\begin{align*}
\bar F \left(\left(\frac{m}{\pi}+\beta\right) k\log k\right)&\ge F_1\left(\left(\frac{m}{\pi}+\beta\right) k\log k\right)>0,\\
\bar F \left(\left(\frac{m}{\pi}-\beta\right)k\log k\right)&\le F_2\left(\left(\frac{m}{\pi}-\beta\right)k\log k\right)<0
\end{align*}
and
\[
\bar F \left(\left(\frac{m}{\pi}+\frac{\beta}{2}\right) k\log k\right)-\bar F \left(\left(\frac{m}{\pi}+\beta\right) k\log k\right)
\ge \frac{C}{k^{2m}\log^{2m}k}>0.
\]
Hence $\bar F$ possesses a critical point (a maximum point)
\[
\bar{r}=\left(\frac{m}{\pi}+o(1)\right)k\log k,
\] 
in the interior of $\S_k$.
\\
Finally it is easy to check that \eqref{eq:max} is achieved by some $r_k$, which is in the interior of $S_k$, and so we infer that $r_k$ is a critical point of $F(r)$. As a consequence, we can conclude that $z_{r_k}+w(r_k)$ is a solution of \eqref{eq:SP'}. This prove the existence of infinitely many non-trivial non radial solutions of \eqref{eq:SP'}. 
\\
In order to get positive solutions, it suffices to repeat the whole procedure for the functional  
$$
I_{+}(u)=\frac{1}{2}\|u\|^2+\frac{1}{4}\int_{\R^3}V(x)\phi_u u^2\, dx-\frac{1}{p+1}\int_{\R^3}|u^+|^{p+1}\, dx,
$$
where $u^+=\max\{0,u\}$. 
It can be checked that Proposition \ref{ausiliaria}, Lemmas \ref{derivprimaoper}, \ref{inv} and Proposition \ref{Izr} can be applied to $I_+$. Therefore we get infinitely many non-radial and non-negative solutions for the problem
$$-\Delta u+u+V(x)\phi_u u= (u^+)^p.$$ Keeping in mind that $\phi_u>0$ when $u\neq 0$ then the maximum principle allows to find infinitely many positive solutions of \eqref{eq:SP'} of the form $z_{r_k}+w(r_k)$, and so the proof of Theorem \ref{principale2} is concluded.

\appendix

\section{Appendix}
\def\theequation{A.\arabic{equation}}\makeatother
\setcounter{equation}{0}

In this section we prove the following result.
\begin{prop}\label{Izr}
For a small $\sigma>0$ there holds
\begin{equation*}
I(z_r)=k\left[C_0+\frac{B_1}{r^{2m}}+ \frac{B_2 k\log k}{r^{2m+1}}-B_3 e^{-\frac{2\pi r}{k}}+O\left(\frac{1}{k^{2m+\sigma}}\right)\right],
\end{equation*}
where $C_0, B_1, B_2, B_3$ are positive constants.
\end{prop}
First we give some preliminary.\\
We recall that
\begin{equation}\label{eq:Iz}
I(z_r)=\frac{1}{2}\int_{\mathbb R^3}\left[|\nabla z_r|^2+z_r^2\right]dx+\frac{1}{4}\int_{\mathbb R^3}V(x)\phi_{z_r}z_r^2dx-\frac{1}{p+1}\int_{\mathbb R^3}|z_r|^{p+1}dx.
\end{equation}
By making simple computations we find
\begin{align}
\int_{\mathbb R^3}V(x)\phi_{z_r}z_r^2dx
&= \sum_{i=1}^k \int_{\mathbb R^3}V(x)\phi_{U_{P_i}}z_r^2dx
+\sum_{i\neq j}\int_{\mathbb R^3}V(x)\phi_{z_r}U_{P_i}U_{P_j}dx\nonumber
\\
&=\sum_{i=1}^k\int_{\mathbb R^3}V(x)\phi_{U_{P_i}}U_{P_i}^2dx
+\sum_{i\neq j}\int_{\mathbb R^3}V(x)\phi_{U_{P_i}}U_{P_j}^2dx\nonumber
\\
&\quad+\sum_{i=1}^k\sum_{j\neq l}\int_{\mathbb R^3}V(x)\phi_{U_{P_i}}U_{P_j}U_{P_l}dx
+\sum_{i\neq j}\int_{\mathbb R^3}V(x)\phi_{z_r}U_{P_i}U_{P_j}dx\nonumber 
\\
&=k\int_{\mathbb R^3}V(x)\phi_{U_{P_1}}U_{P_1}^2dx
+\sum_{i\neq j}\int_{\mathbb R^3}V(x)\phi_{U_{P_i}}U_{P_j}^2dx\nonumber 
\\
&\quad+\sum_{i=1}^k\sum_{j\neq l}\int_{\mathbb R^3}V(x)\phi_{U_{P_i}}U_{P_j}U_{P_l}dx
+\sum_{i\neq j}\int_{\mathbb R^3}V(x)\phi_{z_r}U_{P_i}U_{P_j}dx. \label{eq:Vzz}
\end{align}

\begin{lem}\label{A2}
There holds:
\[
\int_{\mathbb R^3}V(x)\phi_{U_{P_1}}U_{P_1}^2dx
=\frac{B_1}{r^{2m}}+O\left(\frac{1}{k^{2m+ \t}\log^{2m+ \t}k}\right),
\]
where $B_1=\int_{\mathbb R^3}\int_{\mathbb R^3}\frac{U^2(x)U^2(y)}{|x-y|}dxdy$.
\end{lem}

\begin{proof}
Let $\a\in (0,1)$, we have 
\begin{align*}
\int_{\mathbb R^3}V(x)\phi_{U_{P_1}}U_{P_1}^2dx
&=\int_{\mathbb R^3}V(x+P_1)\phi_{U_{P_1}}(x+P_1)U^2(x)dx
\\
&=\int_{B_{\a r}}V(x+P_1)\phi_{U_{P_1}}(x+P_1)U^2(x)dx 
+O(e^{-(1-\tau)r})
\\
&=\int_{B_{\alpha r}}\int_{\mathbb R^3}\frac{V(x+P_1)V(y+P_1)}{|x-y|}U^2(x)U^2(y)dxdy
\\
&\quad+O(e^{-(1-\tau)r})
\\
&=\int_{B_{\alpha r}}\int_{B_{\a r}}\frac{V(x+P_1)V(y+P_1)}{|x-y|}U^2(x)U^2(y)dxdy
\\
&\quad+O(e^{-(1-\tau)r}).
\end{align*}
By \eqref{eq:vp} and \eqref{eq:zp}, since $|P_1|=r$, we have:
\begin{align*}
&\hspace{-2cm}\int_{B_{\a r}}\int_{B_{\a r}}\frac{V(x+P_1)V(y+P_1)}{|x-y|}U^2(x)U^2(y)dxdy
\\
&=\frac{1}{r^{2m}}\int_{B_{\a r}}\int_{B_{\a r}}\frac{U^2(x)U^2(y)}{|x-y|}dxdy
+O\left(\frac{1}{r^{2m+\t}}\right)
\\
&=\frac{B_1}{r^{2m}}+O\left(\frac{1}{k^{2m+\t}\log^{2m+\t}k}\right).
\end{align*}
\end{proof}

\begin{lem}\label{A3}
For a suitable $\s>0$, we have
\begin{equation}\label{pocoimportanti1}
\sum_{i=1}^k\sum_{j\neq l}\int_{\mathbb R^3}V(x)\phi_{U_{P_i}}U_{P_j}U_{P_l}dx=k\cdot O\left(\frac{1}{k^{2m+\sigma}}\right).
\end{equation}
\end{lem}

\begin{proof}
By using H\"{o}lder inequality we obtain
\begin{equation*}
\sum_{i=1}^k \sum_{j\neq l}\int_{\mathbb R^3}V(x)\phi_{U_{P_i}}U_{P_j}U_{P_l}\,dx\le C\sum_{i=1}^k \sum_{j\neq l}\|\phi_{U_{P_i}}\|_{D^{1, 2}}\left(\int_{\mathbb R^3}\left(U_{P_j}U_{P_l}\right)^{\frac65}\,dx\right)^{\frac56}.
\end{equation*}
Then, by Lemma \ref{berlions}, 
\beq\label{primopezzo}
\sum_{i=1}^k \sum_{j\neq l}\int_{\mathbb R^3}V(x)\phi_{U_{P_i}}U_{P_j}U_{P_l}\,dx
\le k\cdot C \|\phi_{U_{P_1}}\|_{D^{1, 2}}\sum_{j\neq l}e^{-(1-\delta)|P_j-P_l|}.
\eeq
Let us evaluate $\|\phi_{U_{P_1}}\|_{D^{1, 2}}$.
\begin{align*}
\|\phi_{U_{P_1}}\|_{D^{1, 2}}^2
&=\int_{\mathbb R^3}|\nabla \phi_{U_{P_1}}|^2\, dx=\int_{\mathbb R^3}V(x)\phi_{U_{P_1}}U_{P_1}^2\,dx
\\
&\le  C\|\phi_{U_{P_1}}\|_{D^{1, 2}}\left(\int_{\mathbb R^3}\left(V(x)U_{P_{1}}^2\right)^{6/5}\, dx\right)^{5/6}
\\
&\le  C\|\phi_{U_{P_1}}\|_{D^{1, 2}}\left(\int_{\mathbb R^3}\left(V(x+P_1)U^2\right)^{6/5}\, dx\right)^{5/6}
\\
&\le  C\|\phi_{U_{P_1}}\|_{D^{1, 2}}\left(\int_{B_{\alpha r}}\left(V(x+P_1)U^2\right)^{6/5}\, dx
+O(e^{-(1-\tau)r})\right)^{5/6},
\end{align*}
where $0<\alpha<1$ and $\tau>0$ sufficiently small.
\\
By (\ref{dis}), we have
\begin{equation*}
\|\phi_{U_{P_1}}\|_{D^{1, 2}}\le C\left(\underbrace{\int_{B_{\alpha r}}\left(V(x+P_1)U^2\right)^{6/5}\, dx}_{(I)}\right)^{5/6}
+O(e^{-\frac 56(1-\tau)r}).
\end{equation*}
Using again (\ref{dis}) we find
\begin{align*}
(I)
&= \int_{B_{\alpha r}}\left[\frac{1}{|x+P_1|^{m}}+O\left(\frac{1}{|x+P_1|^{m+\theta}}\right)\right]^{6/5}U^{12/5}\, dx
\\
&\le  \int_{B_{\alpha r}}\frac{1}{|x+P_1|^{\frac{6m}{5}}}U^{12/5}\, dx
+O\left(\int_{B_{\alpha r}}\frac{1}{|x+P_1|^{\frac{6}{5}(m+\theta)}}U^{12/5}\, dx\right)
\\
&=\frac{1}{|P_1|^{\frac{6}{5}m}}\int_{B_{\alpha r}}U^{12/5}\, dx+O\left(\frac{1}{|P_1|^{\frac{6}{5}(m+\t)}}\right)
\\
&=\frac{1}{r^{\frac{6}{5}m}}\int_{\mathbb R^3}U^{12/5}\, dx+O\left(\frac{1}{r^{\frac{6}{5}(m+\t)}}\right).
\end{align*}
\\
Hence
\begin{equation*}
\|\phi_{U_{P_1}}\|_{D^{1, 2}}\le C_1\frac{1}{r^{m}}+O\left(\frac{1}{r^{m+\t}}+e^{-\frac 56(1-\tau)r}\right)\le C_1\frac{1}{r^{m}}+O\left(\frac{1}{r^{m+\t}}\right).
\end{equation*}
Then from (\ref{primopezzo}) we obtain
\begin{equation*}
\sum_{i=1}^k \sum_{j\neq l}\int_{\mathbb R^3}V(x)\phi_{U_{P_i}}U_{P_j}U_{P_l}\,dx\le  k\cdot C\left(\frac{C_2}{r^m}+O\left(\frac{1}{r^{m+\t}}\right)\right)\sum_{j\neq l}e^{-(1-\delta)|P_j-P_l|}.
\end{equation*}
Since
\begin{equation}\label{epipl}
\sum_{j\neq l}e^{-(1-\delta)|P_j-P_l|}=k\cdot \sum_{j=2}^{k}e^{-(1-\delta)|P_1-P_j|}
\end{equation}
and 
\begin{equation}\label{p1pj}
\sum_{j=2}^{k}e^{-(1-\delta)|P_1-P_j|}\le C e^{-(1-\delta)\frac{2r\pi}{k}},
\end{equation}
recalling that $r\in S_k$, we find
\begin{align*}
\sum_{i=1}^k\sum_{j\neq l}\int_{\mathbb R^3}V(x)\phi_{U_{P_i}}U_{P_j}U_{P_l}\, dx 
&=k \cdot O\left(\frac{1}{k^{m-1}\log^{m}k}\cdot k^{-(1-\delta)(2m-2\pi\beta)}\right)
\\
&=k\cdot O\left(\frac{1}{k^{2m+\sigma}}\right),
\end{align*}
with $$\sigma:=\delta(2\pi\beta-2m)+m(1-\eta)-1-2\pi\beta,$$ where $\eta>0$ small is such that $\log k \ge k^{-\eta}$ for $k$ large. Therefore, since $m>\frac{3}{2}$, $\sigma>0$ if $\delta$, $\beta$ and $\eta$ are sufficiently small. 
\end{proof}

\begin{lem}\label{A4}
For a suitable $\s>0$, we have
\begin{equation*}
\sum_{i\neq j}\int_{\mathbb R^3}V(x)\phi_{z_r}U_{P_i}U_{P_j}dx=k\cdot O\left(\frac{1}{k^{2m+\sigma}}\right).
\end{equation*}
\end{lem}

\begin{proof}
We compute
\begin{align*}
\sum_{i\neq j}\int_{\mathbb R^3}V(x)\phi_{z_r}U_{P_i}U_{P_j}\, dx
&= \underbrace{\sum_{l=1}^k\sum_{i\neq j}\int_{\mathbb R^3}V(x)\phi_{U_{P_l}}U_{P_i}U_{P_j}\, dx}_{(A_1)}
\\
&\quad+\underbrace{\sum_{l\neq t}\sum_{i\neq j}\int_{\mathbb R^3}V(x)\phi_{l, t}U_{P_i}U_{P_j}\, dx}_{(A_2)},
\end{align*}
where $\phi_{l,t}\in D^{1,2}(\mathbb{R}^3)$ is the unique solution of 
$-\Delta \phi=V(x)U_{P_l}U_{P_t}.$
By (\ref{pocoimportanti1}) it follows immediately that $$(A_1)=k\cdot O\left(\frac{1}{k^{2m+\sigma}}\right).$$
Now, using Lemma \ref{berlions}, (\ref{epipl}) and (\ref{p1pj}), since $r\in S_k$, we find
\begin{align*}
(A_2)&\le  C\sum_{l\neq t}\|\phi_{l, t}\|_{D^{1, 2}}
\cdot \sum_{i\neq j}\left(\int_{\mathbb R^3}\left(U_{P_i}U_{P_j}\right)^{6/5}\, dx\right)^{5/6}
\\
&\le  C_1 \sum_{l\neq t}\left(\int_{\mathbb R^3}\left(U_{P_l}U_{P_t}\right)^{6/5}\, dx\right)^{5/6}
\cdot \sum_{i\neq j}\left(\int_{\mathbb R^3}\left(U_{P_i}U_{P_j}\right)^{6/5}\, dx\right)^{5/6}
\\
&\le \sum_{l\neq t}e^{-(1-\delta)|P_l-P_t|}
\cdot \sum_{i\neq j}e^{-(1-\delta)|P_i-P_j|}
\\
&\le k\cdot k\cdot e^{-(1-\delta)\frac{4r\pi}{k}}
\\
&=  k\cdot O\left(\frac{1}{k^{2m+\sigma}}\right),
\end{align*}
where
$$\sigma:=\delta(4\beta\pi-4m)+2m-4\beta\pi-1.$$
Hence, since $m>1$, $\sigma>0$ if $\delta>0$ and $\beta>0$ are sufficiently small. 
\end{proof}

\begin{lem}\label{A5}
For a suitable $\s>0$, we have
\[
\sum_{i\neq j}\int_{\mathbb R^3}\!\!V(x)\phi_{U_{P_i}}U_{P_j}^2\, dx=
k\left[\frac{C}{r^{2m}}\sum_{j=2}^k\frac{1}{|P_1-P_j|}
+\frac{1}{r^{2m+\s}}O\left(\sum_{j=2}^k\frac{1}{|P_1-P_j|}\!\right)\!\right].
\]
\end{lem}

\begin{proof}
First of all, let us observe that
\beq\label{eq:ij1j}
\sum_{i\neq j}\int_{\mathbb R^3}V(x)\phi_{U_{P_i}}U_{P_j}^2\, dx
=k \sum_{j=2}^k\int_{\mathbb R^3}V(x)\phi_{U_{P_1}}U_{P_j}^2\, dx.
\eeq
If $0<\a<1$, we have
\begin{align*}
\int_{\mathbb R^3}V(x)\phi_{U_{P_1}}U_{P_j}^2dx
&=\int_{\mathbb R^3}V(x+P_j)\phi_{U_{P_1}}(x+P_j)U^2(x)dx
\\
&=\int_{B_{\a r}}V(x+P_j)\phi_{U_{P_1}}(x+P_j)U^2(x)dx 
+O(e^{-(1-\tau)r})
\\
&=\int_{B_{\alpha r}}\int_{\mathbb R^3}\frac{V(x+P_j)V(y+P_1)}{|x-y+P_j-P_1|}U^2(x)U^2(y)dxdy
\\
&\quad+O(e^{-(1-\tau)r})
\\
&=\int_{B_{\alpha r}}\int_{B_{\a r}}\frac{V(x+P_j)V(y+P_1)}{|x-y+P_j-P_1|}U^2(x)U^2(y)dxdy
\\
&\quad+O(e^{-(1-\tau)r}).
\end{align*}
We claim that 
\beq\label{eq:xuu}
\int_{B_{\a r}}\int_{B_{\a r}}\frac{|x|U^2(x)U^2(y)}{|x-y+P_j-P_1|}dxdy=O\left(\frac{1}{|P_1-P_j|}\right).
\eeq
Indeed, as in \cite[Lemma 3.2]{DW},
since $$\irt \frac{|x||y-P_j+P_1|U^2(x)}{|x-y+P_j-P_1|}dx \le C,$$
 we have
\begin{align*}
\int_{B_{\a r}}\int_{B_{\a r}}\frac{|x|U^2(x)U^2(y)}{|x-y+P_j-P_1|}dxdy
&\le\irt U^2(y) \left(\irt \frac{|x|U^2(x)}{|x-y+P_j-P_1|}dx \right)dy
\\
&\le C_1 \irt \frac{U^2(y)}{|y-P_j+P_1|}dx
\\
&\le C_2\frac{1}{|P_j-P_1|}=O\left(\frac{1}{|P_1-P_j|}\right).
\end{align*}

Analogously, we can prove that 
\begin{align}
\int_{B_{\a r}}\int_{B_{\a r}}\frac{|y|U^2(x)U^2(y)}{|x-y+P_j-P_1|}dxdy=&O\left(\frac{1}{|P_1-P_j|}\right),\label{eq:yuu}
\\
\int_{B_{\a r}}\int_{B_{\a r}}\frac{|x||y|U^2(x)U^2(y)}{|x-y+P_j-P_1|}dxdy=&O\left(\frac{1}{|P_1-P_j|}\right).\label{eq:xyuu}
\end{align}
Therefore, since $|P_1|=|P_j|=r$, by \eqref{eq:vp} and \eqref{eq:zp}, together with \eqref{eq:xuu}, \eqref{eq:yuu} and \eqref{eq:xyuu}, we have:
\begin{align*}
\int_{\mathbb R^3}V(x)\phi_{U_{P_1}}U_{P_j}^2\, dx
&=\frac{1}{r^{2m}}\irt\irt\frac{U^2(x)U^2(y)}{|x-y+P_j-P_1|}dxdy
\\
&\quad+\frac{1}{r^{2m+\s}}O\left(\frac{1}{|P_1-P_j|}\right).
\end{align*}
Hence, by \cite{DW}, we infer that
\[
\int_{\mathbb R^3}V(x)\phi_{U_{P_1}}U_{P_j}^2\, dx
=\frac{C}{r^{2m}}\frac{1}{|P_1-P_j|}
+\frac{1}{r^{2m+\s}}O\left(\frac{1}{|P_1-P_j|}\right),
\]
and by \eqref{eq:ij1j}, we get the conclusion.
\end{proof}
\begin{lem}\label{A6}
There holds:
\begin{equation*}
\sum_{i=2}^k \frac{1}{|P_1-P_i|}=\frac{k}{\pi r}\log k+o(1).
\end{equation*}
\end{lem}

\begin{proof}
First of all, let us observe 
\begin{equation*}
\sum_{i=2}^k\frac{1}{|P_1-P_i|}
=\frac{1}{2r}\sum_{i=2}^k \frac{1}{\sin\frac{(i-1)\pi}{k}}
=\frac{1}{r}\left(\sum_{i=1}^{\left[\frac{k-1}{2}\right]}\frac{1}{\sin\frac{i\pi}{k}}+ \frac{1+(-1)^k}{4}\right).
\end{equation*}
Then it is sufficient to prove that
\begin{equation}
\lim_k \left(\frac{1}{k\log k}
\sum_{i=1}^{\left[\frac{k-1}{2}\right]}{\frac{1}{\sin \frac{i\pi}{k}}}\right)
=\frac{1}{\pi}
\label{eq:lim}
\end{equation}
Since, for every $i=1,\ldots,\left[\frac{k-1}{2}\right]-1$ and $s\in[i,i+1]$,
\[
 	\sin \frac{i\pi}{k}\le  \sin \frac{s\pi}{k} \le \sin \frac{(i+1)\pi}{k},
\]
%
then
\[
 	\frac{1}{\sin \frac{(i+1)\pi}{k}}=\int_i^{i+1}\frac{1}{\sin \frac{(i+1)\pi}{k}}ds
 	\le \int_i^{i+1}\frac{1}{\sin \frac{s\pi}{k}}ds
 	\le \int_i^{i+1}\frac{1}{\sin \frac{i\pi}{k}}ds=\frac{1}{\sin \frac{i\pi}{k}}.
\]
Hence, adding on $i$, we have that
\[
 	\frac{1}{\sin \left(\left[\frac{k-1}{2}\right]\frac{\pi}{k}\right)}
 	\le
 	\sum_{i=1}^{\left[\frac{k-1}{2}\right]}\frac{1}{\sin \frac{i\pi}{k}}
 	-\int_1^{\left[\frac{k-1}{2}\right]}\frac{1}{\sin \frac{s\pi}{k}}ds
 	\le
 	\frac{1}{\sin \frac{\pi}{k}}
\]
and so
\begin{equation*}
\lim_k \left[\frac{1}{k\log k}\left(
 	\sum_{i=1}^{\left[\frac{k-1}{2}\right]}\frac{1}{\sin \frac{i\pi}{k}}
 	-\int_1^{\left[\frac{k-1}{2}\right]}\frac{1}{\sin \frac{s\pi}{k}}ds\right)\right]=0.
\label{eq:diff0}
\end{equation*}
Therefore \eqref{eq:lim} follows, being
\[
\lim_k \left(\frac{1}{k\log k}\int_1^{\left[\frac{k-1}{2}\right]}\frac{1}{\sin \frac{s\pi}{k}}ds\right)
=\frac{1}{\pi}\lim_k \left(\frac{1}{\log k}\log\left(\frac{\tan \left(\left[\frac{k-1}{2}\right]\frac{\pi}{2k}\right)}{\tan \frac{\pi}{2k}}\right)\right)
=\frac{1}{\pi}.
\]
\end{proof}

We are now ready to prove Proposition \ref{Izr}.

\begin{proof}[Proof of Proposition \ref{Izr}]
By \eqref{eq:Vzz} and using Lemmas \ref{A2}, \ref{A3}, \ref{A4}, \ref{A5}, \ref{A6} we find
\begin{align*}
\frac{1}{4}\int_{\mathbb R^3}V(x)\phi_{z_r}z_r^2\, dx&= k\left[\frac{B_1}{r^{2m}}+\frac{B_2}{r^{2m}}\sum_{i=2}^k \frac{1}{|P_1-P_i|}\right]\\
&\quad +k\left[\frac{1}{r^{2m+\t}}O\left(\sum_{i=2}^k\frac{1}{|P_1-P_i|}\right)+O\left(\frac{1}{k^{2m+\sigma}}\right)\right]\\
&= k\left[\frac{B_1}{r^{2m}}+\frac{B_2}{r^{2m}}\left(\frac{k\log k}{r}+o\left(1\right)\right)
+O\left(\frac{1}{k^{2m+\sigma}}\right)\right]
\\
&= k\left[\frac{B_1}{r^{2m}}+\frac{B_2 k \log k}{r^{2m+1}}+O\left(\frac{1}{k^{2m+\sigma}}\right)\right].
\end{align*}
Then, by \eqref{eq:Iz}, we have
\begin{align*}
I(z_r)&=\frac12\sum_{j=1}^k\sum_{i=1}^k\int_{\mathbb R^3}U_{P_j}^p U_{P_i}\, dx+k\left[\frac{B_1}{r^{2m}}+\frac{B_2k\log k}{r^{2m+1}}+O\left(\frac{1}{k^{2m+\sigma}}\right)\right]\\
&\quad-\frac{1}{p+1}\int_{\mathbb R^3}|z_r|^{p+1}\, dx\\
&=\frac k2\int_{\mathbb R^3}U^{p+1}\, dx+\frac k2\sum_{i=2}^k\int_{\mathbb R^3}U_{P_1}^p U_{P_i}\, dx+k\left[\frac{B_1}{r^{2m}}+\frac{B_2k\log k}{r^{2m+1}}+O\left(\frac{1}{k^{2m+\sigma}}\right)\right]\\
&\quad-\frac{1}{p+1}\int_{\mathbb R^3}|z_r|^{p+1}\, dx.
\end{align*}
As done in \cite[Proof of Proposition A.3]{WeiYan}, for all $r\in S_k$, one can prove 
\begin{equation*}
\int_{\mathbb R^3}|z_r|^{p+1}\, dx=k \left(\int_{\mathbb R^3}U^{p+1}\, dx+(p+1)\sum_{i=2}^k\int_{\mathbb R^3}U_{P_1}^p U_{P_i}\, dx+O\left(\frac{1}{k^{2m+\sigma}}\right)\right).
\end{equation*}
At the end we find 
\begin{align*}
I(z_r)&=k\left[\left(\frac{1}{2}-\frac{1}{p+1}\right)\int_{\mathbb R^3}U^{p+1}\, dx+\frac{B_1}{r^{2m}}+\frac{B_2k\log k}{r^{2m+1}}\right.\\
&\qquad\left.-B_3\sum_{i=2}^k\int_{\mathbb R^3}U_{P_1}^p U_{P_i}\, dx+O\left(\frac{1}{k^{2m+\sigma}}\right)\right]\\
&=k\left[C_0+\frac{B_1}{r^{2m}}+\frac{B_2k\log k}{r^{2m+1}}-B_3\sum_{i=2}^k\int_{\mathbb R^3}U_{P_1}^p U_{P_i}\, dx+O\left(\frac{1}{k^{2m+\sigma}}\right)\right].
\end{align*}

\end{proof}


\begin{thebibliography}{99}
\bibitem{ACR}
A. Ambrosetti, E. Colorado, D. Ruiz, {\it Multi-bump solitons to linearly coupled systems of nonlinear Schr\"odinger equations},  Calc. Var. Partial Differential Equations, {\bf 30}, (2007), 85--112.

\bibitem{AmbrMalc}
A. Ambrosetti, A. Malchiodi, Perturbation Methods and Semilinear Elliptic Problems on ${\mathbb R}^{n}$,
Birkh\"{a}user Verlag, 2005.

\bibitem{A}
A. Azzollini, {\it Concentration and compactness in nonlinear Schr\"odinger-Poisson system}, preprint.

\bibitem{ADP}
A. Azzollini, P. d'Avenia, A. Pomponio, {\it On the Schr\"odinger-Maxwell equations under the effect of a general nonlinear term}, Ann. Inst. H. Poincar\'e Anal. Non Lin\'eaire, {\bf 27}, (2010), 779--791.

\bibitem{AP}
A. Azzollini, A. Pomponio, {\it Ground state solutions for the nonlinear Schr\"{o}dinger-Maxwell equations}, J. Math. Anal. Appl., {\bf 345}, (2008), 90--108.

\bibitem{BenciFortunato1}
V. Benci, D. Fortunato, {\it An eigenvalue problem for the Schr\"{o}dinger--Maxwell equations}, 
Topol. Methods Nonlinear Anal., {\bf 11}, (1998), 283--293.



\bibitem{CeramiVaira}
G. Cerami, G. Vaira, {\it Positive solutions for some non autonomous Schr\"{o}dinger-Poisson systems}, 
J. Differential Equations, {\bf 248},(2010), 521--543.

\bibitem{DM1}
T. D'Aprile, D. Mugnai, {\it Solitary waves for nonlinear
Klein-Gordon-Maxwell and Schr\"{o}dinger-Maxwell equations},
Proc. Roy. Soc. Edinburgh Sect. A, {\bf 134}, (2004), 893--906.




\bibitem{DW}
T. D'Aprile, J. Wei, {\it Standing waves in the Maxwell-Schr\"odinger equation and an optimal configuration problem},  Calc. Var. Partial Differential Equations, {\bf 25}, (2006), 105--137.


\bibitem{DA}
P. d'Avenia, {\it Non-radially symmetric solutions of nonlinear Schr\"odinger equation coupled with Maxwell equations}, Adv. Nonlinear Stud., {\bf 2}, (2002),  177--192.


\bibitem{KW}
X, Kang, J. Wei, {\it On interacting bumps of semi-classical states of nonlinear Schr\"odinger equations}, 
Adv. Differential Equations, {\bf 5}, (2000),  899--928.


\bibitem{K1}
H. Kikuchi, {\it On the existence of a solution for elliptic system related to the Maxwell-Schr\"odinger equations},
Nonlinear Anal., Theory Methods Appl., {\bf 67}, (2007), 1445--1456.


\bibitem{K2}
H. Kikuchi, {\it Existence and stability of standing waves for Schr\"{o}dinger-Poisson-Slater equation},
Adv. Nonlinear Stud., {\bf 7}, (2007), 403--437.


\bibitem{Kwong}
M. K. Kwong, {\it Uniqueness of positive solutions of $-\Delta u+u=u^p$ in $\mathbb R^N$}, 
Arch. Rational Mech. Anal., {\bf 105},  (1989), 243--266.

\bibitem{IanniVaira}
I. Ianni, G. Vaira, {\it On Concentration of Positive Bound States
for the Schr\"{o}dinger-Poisson Problem with Potentials}, 
Adv. Nonlinear Stud., {\bf 8}, (2008), 573--595.


\bibitem{Ru}
D. Ruiz, \textit{The Schr\"odinger-Poisson equation under the
effect of a nonlinear local term}, Journ. Func. Anal., {\bf 237},
(2006), 655--674.



\bibitem{RuizVaira}
D. Ruiz, G. Vaira, {\it Cluster solutions for the Schr{\"o}dinger-Poisson-Slater problem around a local minimum of the potential}, to appear on Rev. Mat. Iberoamericana.

\bibitem{WWY}
L. Wang, J. Wei, S. Yan  {\it A Neumann Problem with Critical Exponent in Non-convex Domains and Lin-Ni's Conjecture}, 
to appear on Trans. Amer. Math. Soc. 


\bibitem{WZ}
Z. Wang, H.S. Zhou,
{\it Positive solution for a nonlinear stationary Schr\"odinger-Poisson system in $\RT$},
Discrete Contin. Dyn. Syst., {\bf 18}, (2007), 809--816.



\bibitem{WeiYan}
J. Wei, S. Yan, {\it Infinitely many positive solutions for the nonlinear Schr\"{o}dinger equations in $\mathbb R^N$}, Calc. Var. Partial Differential Equations, {\bf 37}, (2010), 423--439.


\bibitem{WY2}
J. Wei, S. Yan, {\it Infinitely many solutions for the prescribed scalar curvature problem on $\mathbb{S}^N$}, Journ. Func. Anal., {\bf 258}, (2010), 3048--3081. 


\bibitem{ZZ}
L. Zhao, F. Zhao,
{\it On the existence of solutions for the Schr\"odinger-Poisson equations},
J. Math. Anal. Appl., {\bf 346}, (2008), 155--169.


\end{thebibliography}
\end{document}